\documentclass[11pt, a4paper]{article}
\usepackage{latexsym}
\usepackage{indentfirst}
\usepackage{graphicx}
\usepackage{amsmath}
\usepackage{amssymb}

\newtheorem{lem}{Lemma}[section]

\newtheorem{thm}{Theorem}[section]

\begin{document}
\title{On some extremal problems in spaces of harmonic functions}
\author{Milo\v s Arsenovi\'c$\dagger$, Romi F. Shamoyan  \\
University of Belgrade, Belgrade, Serbia\\
Bryansk University, Bryansk, Russia \\
arsenovic@matf.bg.ac.rs, rshamoyan@gmail}

\date{}
\maketitle

\begin{abstract}
We give solutions to some extremal problems involving distance function in mixed normed spaces of
harmonic functions on the unit ball in $\mathbb R^n$.
\end{abstract}

\thanks{$\dagger$ Supported by Ministry of Science, Serbia, project OI174017}

\section{Introduction}
The classical duality approach to distance problems in spaces of analytic functions is based on a
functional analysis result:
$${\rm dist}_X (f, Y) = \inf_{g \in Y} \| f - g \|_X= \sup_{\phi \in Y^\perp, \| \phi \| \leq 1} | \phi (f) |,$$
where $Y^\perp \subset X^\star$ is the annihilator of a subspace $Y$ of a normed space $X$. This method was employed
in the case of Hardy spaces, \cite{Kh1}, \cite{Du}. More recently, extremal problems in the Bergman spaces were considered in \cite{KS}. The aim of this paper is to obtain extremal distance results in spaces of harmonic
functions on the unit ball $\mathbb B \subset \mathbb R^n$ using direct methods introduced by R. Zhao, \cite{RZ}, \cite{WX}. We note that spaces of harmonic functions on the unit ball $\mathbb B$ were extensively studied, see \cite{DS}, \cite{GKU} and \cite{Ka} and references therein.

In this paper letter $C$ designates a positive constant which can change its value even in the same chain of inequalities. Given real expressions $A$ and $B$, we write $A\lesssim B$ if there is a constant $C \geq 0$ such that $A \leq CB$. If $A \lesssim B$ and $B \lesssim A$ then we write $A \asymp B$. The integer part of a real number $\alpha$ is denoted by $[\alpha]$, and the fractional part of $\alpha$ is denoted by $\{\alpha \} = \alpha -[\alpha]$.

We use standard notation: $\mathbb S = \partial \mathbb B$ is the unit sphere in $\mathbb R^n$, for $x \in
\mathbb R^n$ we have $x = rx'$, where $r = |x| = \sqrt{\sum_{j=1}^n x_j^2}$ and $x' \in \mathbb S$. Normalized
Lebesgue measure on $\mathbb B$ is denoted by $dx = dx_1 \ldots dx_n = r^{n-1}dr dx'$ so that $\int_{\mathbb B} dx = 1$. The space of all harmonic functions on $\mathbb B$ is denoted by $h(\mathbb B)$.

We consider harmonic weighted Bergman spaces $A^p_\alpha (\mathbb B)$ on $\mathbb B$ defined for $0 \leq \alpha < \infty$ and $0 < p \leq \infty$ by
$$A^p_\alpha(\mathbb B) = \left\{ f \in h(\mathbb B) : \| f \|_{p, \alpha} = \left( \int_{\mathbb B}
|f(rx')|^p(1-r)^\alpha r^{n-1}dr dx' \right)^{1/p} < \infty \right\},$$
$$A^\infty_\alpha(\mathbb B) = \left\{ f \in h(\mathbb B) : \|f\|_{\infty, \alpha} = \sup_{x \in \mathbb B}
|f(x)|(1-|x|)^\alpha < \infty \right\}.$$
The spaces $A^p_\alpha = A^p_\alpha (\mathbb B)$ are Banach spaces for $1 \leq p \leq \infty$ and complete metric spaces for $0<p<1$.

We are going to use some results on spherical harmonics from \cite{SW} and $A^p_\alpha$ spaces from \cite{DS}.
Let $Y^{(k)}_j$ be the spherical harmonics of order $k$, $1 \leq j \leq \alpha_k$, on $\mathbb S$. Next,
$$Z_{x'}^{(k)}(y') = \sum_{j=1}^{\alpha_k} Y_j^{(k)}(x') \overline{Y_j^{(k)}(y')},$$
are zonal harmonics of order $k$. Note that the spherical harmonics
$Y^{(k)}_j$, ($k \geq 0$, $1 \leq j \leq \alpha_k$) form an orthonormal basis of $L^2(\mathbb S, dx')$. Every
$f \in h(\mathbb B)$ has an expansion
$$f(x) = f(rx') = \sum_{k=0}^\infty r^k Y^{(k)}(x'),$$
where $Y^{(k)} = \sum_{j=1}^{\alpha_k} c_k^j Y_j^k$. Using this expansion one defines, as in \cite{DS}, fractional derivative of order $t \in \mathbb R$ of $f \in h(\mathbb B)$ by the following formula:
$$\mathcal D^t f(rx') = \sum_{k=0}^\infty r^k \frac{\Gamma(k+t+n/2)}{\Gamma(k+n/2)\Gamma(t+n/2)}Y^{(k)}(x').$$

Let
$$Q_\alpha(x, y) = 2 \sum_{k=0}^\infty \frac{\Gamma(\alpha + 1 + k + n/2)}{\Gamma(\alpha + 1) \Gamma(k + n/2)}
r^k \rho^k Z_{x'}^{(k)}(y'), \qquad x = rx', y = \rho y' \in \mathbb B$$
for $\alpha > -1$, where $\Gamma$ is the Euler's Gamma function.
$Q_\alpha(x, y)$ is harmonic in each of the variables separately, $Q_\alpha(x, y) = \overline{Q_\alpha(y, x)}$ and
$\| Q_\alpha (x, \cdot)\|_{A^p_\alpha} = \| Q_\alpha(\cdot, y) \|_{A^p_\alpha} = 1$, see \cite{DS}. In fact, it is a harmonic Bergman kernel for the space $A^p_\alpha$, we have the following theorem from \cite{DS}, see also
\cite{JP}:
\begin{thm}\label{Thm1}
Let $p \geq 1$ and $\alpha \geq 0$. Then for every $f \in A^p_\alpha$ and $x \in \mathbb B$ we have
$$f(x) = \int_0^1 \int_{\mathbb S} (1-\rho^2)^\alpha Q_\alpha(x, y) f(\rho y') \rho^{n-1}d\rho dy',
\qquad y = \rho y'.$$
\end{thm}
The following lemma from \cite{DS} gives estimates for this kernel.
\begin{lem}\label{Lemma1}
1. Let $\alpha > 0$. Then, for $x = rx', y = \rho y' \in \mathbb B$ we have
$$|Q_\alpha(x, y)| \leq C\frac{(1-r\rho)^{-\{ \alpha \}}}{|\rho x - y'|^{n + [\alpha]}} + \frac{C}{(1-r\rho)^{1+\alpha}}.$$

2. Let $\beta > -1$.
$$\int_{\mathbb S} |Q_\beta(rx', y)| dx' \leq \frac{C}{(1-r\rho)^{1+\beta}}, \qquad |y| = \rho.$$

3. Let $m > n-1$, , $0 \leq r < 1$ and $y'\in \mathbb S$. Then
$$\int_{\mathbb S} \frac{dx'}{|rx' - y'|^m} \leq \frac{C}{(1-r)^{m-n+1}}.$$
\end{lem}

Next we note that there is a variant of Theorem \ref{Thm1} for $0<p<1$, Theorem \ref{Thm2} below, which is weaker in the sense that it does not provide integral representation, but stronger in the sense that allows more general weights. First we introduce the class of admissible weight functions.

Let $S$ be the class of all functions $\omega(t) \geq 0$, $0 < t < 1$ such that
$$m_\omega \leq\frac{\omega(\lambda r)}{\omega(r)} \leq M_\omega, \qquad \lambda \in [q_\omega, 1], 0 < r < 1.$$
for some constants $m_\omega, M_\omega, q_\omega \in (0, 1)$. An example of a function $\omega$ in $S$ is
$$\omega(t) = t^\alpha \left(\log\left(\frac{C}{t}\right)\right)^\beta, \qquad \alpha > -1, \beta > 0.$$
For $\omega \in S$ we set $\alpha_\omega = \log m_\omega / \log q_\omega$ and $\beta_\omega = - \log M_\omega/ \log
q_\omega$. To each $\omega \in S$ and $0 < p \leq 1$ we associate a weighted harmonic Bergman space
$$A_\omega^p(\mathbb B) = \left\{ f \in h(\mathbb B) : \int_0^1 \int_{\mathbb S} |f(rx')|^p \omega(1-r) r^{n-1}
dx'dr < \infty \right\}.$$
A characterization of functions in $A_\omega^p(\mathbb B)$ was obtained in \cite{SZ}, in order to formulate this result
it is convenient to introduce some notation. We recall polar coordinates of $x \in \mathbb R^n$:
$x_1 = r \cos \phi_1$, $x_2 = r \sin \phi_1 \cos \phi_2$,...,$x_{n-1} = r \sin \phi_1 \ldots \sin \phi_{n-2} \cos \phi_{n-1}$, $x_n = r \sin \phi_1 \ldots \sin \phi_{n-1}$ where $r > 0$, $0\leq \phi_i < \pi$ for $1 \leq i \leq n-2$
and $-\pi < \phi_{n-1} \leq \pi$. Given $k \geq 0$, $0 \leq l_i < 2^k$, $1 \leq i < n-2$ and $-2^k \leq l_{n-1} < 2^k$
we define $\Delta_{k;l_1, \ldots, l_{n-1}}$ as the set of all $x \in \mathbb B$ such that the polar coordinates of $x$
satisfy: $1-2^{-k} \leq r < 1 - 2^{-k-1}$, $ \pi l_i/2^k \leq \phi_i < \pi (l_i + 1)/2^k$ for $1 \leq i \leq n-2$ and
$\pi l_{n-1}/2^k < \phi_{n-1} \leq \pi (l_{n-1} + 1)/2^k$. We denote the Poisson kernel for the unit ball by $P(x, y)$,
$x \in \mathbb B$, $y \in \mathbb S$. Since it is harmonic in $x$ we can introduce, for $\alpha > 1$,
$$P_\alpha(x, y) = \mathcal D_x^{\alpha -1}P(x, y).$$

\begin{thm}[\cite{SZ}]\label{Thm2}
Let $\omega \in S$, $0< p \leq 1$ and $\alpha > \frac{\alpha_\omega + n}{p}-n$.
Then a function $f \in h(\mathbb B)$ belongs to $A^p_\omega$ if and only if
$$\int_{\mathbb B} (1-|y|^2)^\alpha P_\alpha(x, y) f(y) d\mu(y) < \infty$$
for every Borel measure $\mu$ on $\mathbb B$ such that
$$\sum_{k=0}^\infty \sum_{l_1 = 0}^{2^k-1} \cdots \sum_{l_{n-1} = -2^k}^{2^k-1}
[\mu(\Delta_{k;l_1, \ldots, l_{n-1}})]^p \omega(|\Delta_{k;l_1, \ldots, l_{n-1}}|^{1/n})
|\Delta_{k;l_1, \ldots, l_{n-1}}|^{1-p} < \infty.$$
\end{thm}

We record here for later use the following well known result.

\begin{lem}\label{Lemma2}
Let $\beta > -1$, $s > -1$, $\gamma > 0$ and $0 < p \leq 1$. Then for every increasing function $G(r)$, $0 \leq r < 1$ we have the following estimate:
\begin{eqnarray*}
\int_0^1 \left(\int_0^1 \frac{G(r)(1-r)^\beta dr}{(1-r\rho)^\gamma}\right)^p (1-\rho)^s d \rho & \leq & C \int_0^1 \int_0^1 \frac{G(r)^p (1-r)^{\beta p + p -1}}{(1-r\rho)^{\gamma p}}\cr
& & (1-\rho)^s d\rho dr,
\end{eqnarray*}
where $C$ depends only on $\beta$, $s$, $\gamma$ and $p$.
\end{lem}

For $f \in h(\mathbb B)$ we define $p$-integral means of $f$ by
$$M_p(f, r) = \left( \int_{\mathbb S} |f(rx')|^p dx' \right)^{1/p}, \qquad 0\leq r < 1$$
for $0<p<\infty$, with the usual modification for $p = \infty$. The function $M_p(f, r)$ is increasing in $0\leq r < 1$
for $p \geq  1$ because $|f|^p$ is subharmonic for $1 \leq p < \infty$. We are going to consider the following spaces
of harmonic functions, for $0<p<\infty$, $1 \leq q < \infty$ and $\alpha > 0$:
$$B^{\infty, q}_\alpha = \left\{ f \in h(\mathbb B) : \| f \|_{B^{\infty, q}_\alpha} = \sup_{0\leq r < 1} 
M_q(f, r)(1-r)^\alpha < \infty \right\},$$
$$B^{p,q}_\alpha = \left\{ f \in h(\mathbb B) : \| f \|_{B^{p,q}_\alpha} = \left( \int_0^1 M^p_q(f, r) (1-r)^{\alpha p -1} dr \right)^{1/p} < \infty \right\},$$
$$B^{p,\infty}_\alpha = \left\{ f \in h(\mathbb B) : \| f \|_{B^{p,\infty}_\alpha} = \left( \int_0^1 M^p_\infty(f, r) (1-r)^{\alpha p -1} dr \right)^{1/p} < \infty \right\}.$$
These spaces have obvious (quasi)-norms, with respect to these (quasi)-norms they are Banach spaces or complete metric
spaces. The following elementary lemma is used below.
\begin{lem}\label{Lemma3}
If $\beta > 0$, the for every increasing function $\phi(r)$, $0 \leq r < 1$ we have
$$\sup_{0 \leq r < 1} \phi(r) (1-r)^\beta \leq \beta \int_0^1 \phi(r) (1-r)^{\beta -1} dr.$$
\end{lem}

\textbf{Proof. } For every $r_0 \in [0, 1)$ we have
\begin{align*}
\int_0^1 \phi(r)(1-r)^{\beta -1} dr \geq & \int_{r_0}^1 \phi(r) (1-r)^{\beta-1} dr \cr
& \geq \phi(r_0) \int_{r_0}^1 (1-r)^{\beta-1}dr = \frac{1}{\beta} \phi(r_0) (1-r_0)^\beta. \quad \Box
\end{align*}
Using this lemma and a fact that $M_q(f, r)$ is increasing for $q \geq 1$ we obtain the following estimate, valid for
$q \geq 1$, $0 < p < \infty$ and $\alpha > 0$:
$$\| f \|^p_{B^{\infty, q}_\alpha} = \sup_{r<1}M_q^p(f, r)(1-r)^{\alpha p} \leq \alpha p \int_0^1 M_q^p(f, r)
(1-r)^{\alpha p -1} dr = \alpha p \| f \|^p_{B^{p,q}_\alpha}.$$
This shows that $B^{p,q}_\alpha$ is continuously embedded in $B^{\infty, q}_\alpha$ for the above range of parameters.
Also, it is elementary that $B^{\infty, 1}_\alpha \subset A^1_\alpha$. These two embedding results allow us to use the integral representation from Theorem \ref{Thm1} to functions $f$ in $B^{\infty, 1}_\alpha \supset B^{p, 1}_\alpha$, $p > 0$. Also, the last inclusion leads to a natural problem of finding ${\rm dist}_{B^{\infty, 1}_\alpha} (f, B^{p, 1}_\alpha)$, which we investigate in the next section.

\section{Main results}
This section is devoted to formulations and proofs of our main results on extremal problems in spaces of harmonic
functions. In order to state our first theorem we introduce the following notation: for $\epsilon > 0$, $\alpha > 0$ and $f \in h(\mathbb B)$ we set
$$L_{\epsilon, \alpha}(f) = \left\{ r \in [0, 1) : (1-r)^\alpha \int_{\mathbb S} |f(rx')| dx' \geq \epsilon
\right\}.$$

\begin{thm}\label{Thm3}
Let $1 \leq p < \infty$, $\alpha > 0$. Set, for $f \in B^{\infty, 1}_\alpha$,
\begin{align*}
s_1 = s_1(f) = & {\rm dist}_{B^{\infty, 1}_\alpha} (f, B^{p, 1}_\alpha), \cr
s_2 = s_2(f) = & \inf \left\{ \epsilon > 0 : \int_0^1 \chi_{L_{\epsilon, \alpha}(f)}(r) (1-r)^{-1} dr < \infty
\right\}.
\end{align*}
Then $s_1 \asymp s_2$, with constants involved depending only on $p$, $\alpha$ and $n$.
\end{thm}

\textbf{Proof. } First we prove that $s_1 \geq s_2$. Assume to the contrary that $s_1 < s_2$. Then there are $\epsilon >
\epsilon_1 > 0$ and $f_1 \in B^{p, 1}_\alpha$ such that $\| f - f_1 \|_{B^{\infty, 1}_\alpha} \leq \epsilon_1$ and
$\int_0^1 \chi_{L_{\epsilon, \alpha}}(r)(1-r)^{-1} dr = \infty$. Hence
\begin{eqnarray*}
(1-r)^\alpha \int_{\mathbb S} |f_1(rx')|dx' & \geq & (1-r)^\alpha \int_{\mathbb S} |f(rx')|dx' \cr
&  & - \sup_{0\leq r < 1} (1-r)^\alpha \int_{\mathbb S} |f(rx') - f_1(rx')| dx' \cr
& \geq & (1-r)^\alpha \int_{\mathbb S} |f(rx')| dx' - \epsilon_1,\cr
\end{eqnarray*}
and therefore we have, for $r \in L_{\epsilon, \alpha}(f)$, the following estimate:
\begin{equation}\label{ee2}
(1-r)^\alpha \int_{\mathbb S} |f_1(rx'|dx' \geq (1-r)^\alpha \int_{\mathbb S} |f(rx')|dx' - \epsilon_1
\geq \epsilon - \epsilon_1.
\end{equation}
However, this leads to a contradiction:
\begin{eqnarray*}
(\epsilon - \epsilon_1)^p \int_0^1 \chi_{L_{\epsilon, \alpha}(f)}(r) (1-r)^{-1} dr & \leq & \int_0^1 M_1^p(f_1, r)
(1-r)^{\alpha p -1} dr \cr
& = & \| f_1 \|_{B^{p,1}_\alpha}^p < \infty.
\end{eqnarray*}

Next we prove $s_1 \leq C s_2$. We fix $f \in B^{\infty, 1}_\alpha$ and $\epsilon > 0$ such that the integral appearing in the definition of $s_2(f)$ is finite. Then, using integral representation from Theorem \ref{Thm1}:
\begin{align}\label{dec}
f(x)  = & \int_0^1 \int_{\mathbb S} (1 - \rho^2)^\alpha Q_\alpha(x, y) f(\rho y') \rho^{n-1} d\rho dy' \cr
 = & \int_{L_{\epsilon, \alpha}(f)} \int_{\mathbb S} (1 - \rho^2)^\alpha Q_\alpha(x, y) f(\rho y') \rho^{n-1} d\rho dy' \cr
 + & \int_{I \setminus L_{\epsilon, \alpha}(f)} \int_{\mathbb S} (1 - \rho^2)^\alpha Q_\alpha(x, y) f(\rho y') \rho^{n-1} d\rho dy' \cr
 = & f_1(x) + f_2(x),
\end{align}
where $I = [0, 1)$. To complete the proof it suffices to prove
\begin{equation}\label{cetiri}
\| f_1 \|_{B^{p,1}_\alpha} \leq C_\epsilon \| f \|_{B^{\infty, 1}_\alpha}
\end{equation}
and
\begin{equation}\label{pet}
\| f_2 \|_{B^{\infty, 1}_\alpha} \leq C\epsilon.
\end{equation}

In proving (\ref{cetiri}) it suffices to consider the case $p=1$, since the general case then follows from the estimate
$\| f_1 \|_{B^{p,1}_\alpha} \leq C \| f_1 \|_{B^{1,1}_\alpha}$ from \cite{OF}. Now we have, using Fubini's theorem
and Lemma \ref{Lemma1}:
\begin{eqnarray*}
\| f_1 \|_{B^{1,1}_\alpha} & = & \int_{\mathbb B} (1-\rho)^{\alpha -1} |f_1(\rho x')| d\rho dx' \cr
& \leq & \int_{\mathbb B} (1-\rho)^{\alpha - 1} \int_{L_{\epsilon, \alpha}(f)} \int_{\mathbb S}
\left( \frac{C}{|1-r\rho|^{1+\alpha}} + \frac{C(1-r\rho)^{-\{\alpha \}}}{|r\rho x' - y'|^{n+[\alpha]}} \right) \cr
& \times & (1-r^2)^\alpha |f(ry')|dr dy' d\rho dx' \cr
& \leq & C \int_{L_{\epsilon, \alpha}(f)} \int_{x' \in \mathbb S} (1-r^2)^\alpha |f(ry')| \int_0^1
\frac{(1-\rho)^{\alpha - 1}}{(1-r\rho)^{1+ \alpha}} d\rho dy' dr \cr
& \leq & C \int_{L_{\epsilon, \alpha}(f)} (1-r)^{\alpha - 1} \int_{y' \in \mathbb S} |f(ry')| dy' dr \cr
& \leq & C \| f_1 \|_{B^{\infty, 1}_\alpha} \int_{L_{\epsilon, \alpha}(f)} (1-r)^{-1} dr \cr
& = & C_\epsilon \| f_1 \|_{B^{\infty, 1}_\alpha}.
\end{eqnarray*}

Now we turn to (\ref{pet}). For every $0 \leq \rho < 1$ we have, using Lemma \ref{Lemma1}, the following estimate:
\begin{eqnarray*}
(1-\rho)^\alpha \int_{\mathbb S} |f_2(\rho x')| dx' & \leq & C(1-\rho)^\alpha \int_{\mathbb S}
\int_{I \setminus L_{\epsilon, \alpha}(f)} \int_{\mathbb S} |f(ry')|(1-r^2)^\alpha \cr
& \times & \left( \frac{(1-\rho r)^{-\{ \alpha \}}}{|r \rho x' - y'|^{n + [\alpha]}} + \frac{1}{(1-r\rho)^{1+\alpha}}
\right) dy' dr dx'
\end{eqnarray*}
Now, using Lemma \ref{Lemma1} again and Fubini's theorem we obtain, taking into account definition of the set
$L = L_{\epsilon, \alpha}(f)$
\begin{eqnarray*}
(1-\rho)^\alpha \int_{\mathbb S} |f_2(\rho x')| dx' & \leq &
C(1- \rho)^\alpha \int_{I\setminus L} \int_{\mathbb S}|f(ry')| (1-r^2)^\alpha \frac{dy' dr}{(1-r\rho)^{1+\alpha}} \cr
& \leq & C\epsilon (1-\rho)^\alpha \int_{I\setminus L} (1-r)^{-\alpha}(1-r^2)^\alpha
\frac{dr}{(1-r\rho)^{1+\alpha}} \cr
& \leq & C\epsilon,
\end{eqnarray*}
which ends the proof of Theorem \ref{Thm3}. $\Box$

Our next results is an analogue of the previous theorem for the case $0 < p \leq 1$.
\begin{thm}\label{Thm4}
Let $0 < p \leq 1$, $\alpha > 0$ and $t > \alpha - 1$. For $f \in B^{\infty, 1}_\alpha$ we set
\begin{equation*}
\hat s_1(f) = {\rm dist}_{B^{\infty, 1}_\alpha} (f, B^{p,1}_\alpha),
\end{equation*}
\begin{equation*}
\hat s_2(f) = \inf \left\{ \epsilon > 0: \int_0^1 \left( \int_0^1 \chi_{L_{\epsilon, \alpha}(f)}(r)
\frac{(1-r)^{t-\alpha}}{(1-r\rho)^{t+1}} dr \right)^p (1-\rho)^{p\alpha -1} d\rho < \infty \right\}.
\end{equation*}
Then $\hat s_1(f) \asymp \hat s_2(f)$, with constants involved depending only on $\alpha$, $n$, $p$ and $t$.
\end{thm}

We note that for $p=1$ we have, using Fubini's theorem, $\hat s_2(f) = s_2(f)$. Therefore, for $p=1$, the above two theorems give the same answer to the distance problem.

\textbf{Proof. } Our argument here is similar to the proof of the previous theorem. We first prove that $\hat s_2(f) \leq
\hat s_1(f)$. Assume that $\hat s_1(f) < \hat s_2(f)$. Then there are $0 < \epsilon_1 < \epsilon$ and
$f_1 \in B^{p,1}_\alpha$ such that
\begin{equation}\label{contr}
\int_0^1 \left( \int_0^1 \chi_{L_{\epsilon, \alpha}(f)}(r) \frac{(1-r)^{t - \alpha}}{(1-r\rho)^{t+1}} dr \right)^p
(1-\rho)^{p\alpha - 1} d \rho = \infty
\end{equation}
and $\| f - f_1 \|_{B^{\infty, 1}_\alpha} \leq \epsilon_1$. We can repeat the argument given at the beginning of the proof of Theorem \ref{Thm3} to arrive at estimate (\ref{ee2}) for $r \in L_{\epsilon, \alpha}(f)$, and this implies
\begin{equation}\label{chi}
\chi_{L_{\epsilon, \alpha}(f)}(r) \leq (\epsilon - \epsilon_1)^{-q} \left( \int_{\mathbb S} |f_1(rx')| dx'
\right)^q(1-r)^{\alpha q}, \qquad 0<q<\infty.
\end{equation}
Now Lemma \ref{Lemma2} can be applied to $G(r) = M(f_1, r)$, $\beta = t$, $s = \alpha p - 1$ and $\gamma = t+1$ to obtain, using (\ref{chi}) with $q=1$,
\begin{align*}
& \int_0^1 \left( \int_0^1 \frac{(1-r)^{t-\alpha} \chi_{L_{\epsilon, \alpha}(f)}(r) dr}{(1-r\rho)^{t+1}} \right)^p
(1-\rho)^{\alpha p -1} d\rho \cr
\leq & C \int_0^1 \int_0^1 \left( \int_{\mathbb S} |f_1(rx')dx' \right)^p
\frac{(1-r)^{(t+1)p - 1}}{(1-r\rho)^{(t+1)p}} (1- \rho)^{\alpha p -1} dr d\rho \cr
\leq & C \| f_1 \|_{B^{p,1}_\alpha}.
\end{align*}
But this gives a contradiction with (\ref{contr}).

Next we prove $\hat s_1(f) \leq C \hat s_2(f)$. Let us fix $f \in B^{\infty, 1}_\alpha$ and $\epsilon > 0$ such that
the integral appearing in the definition of $\hat s_2(f)$ is finite. We use the same decomposition $f = f_1 + f_2$ that appear in (\ref{dec}), and the same reasoning as in the proof of Theorem \ref{Thm3} gives $\| f_2 \|_{B^{\infty, 1}_\alpha} \leq C \epsilon$. Therefore it remains to prove that $f_1 \in B^{p, 1}_\alpha$. The following
chain of inequalities
\begin{align*}
\| f_1 \|_{B^{p,1}_\alpha}^p  = & \int_0^1 (1-\rho)^{\alpha p -1} \Big( \int_{\mathbb S} |f_1(\rho x')| dx'
\Big)^p d\rho\cr
\leq & C \int_0^1 (1-\rho)^{\alpha p -1} \Big( \int_{L_{\epsilon, \alpha}(f)} \int_{\mathbb S} (1-r^2)^t
\int_{\mathbb S} \Big( \frac{1}{|r\rho x' - y'|^{n+t}} +  \cr
& + \frac{(1-r\rho)^{-\{t\}}}{|r\rho x' - y'|^{n+[t]}} \Big) |f(ry')| dr dy' \Big)^p d\rho \cr
& \leq \left( \sup_{0 \leq r < 1} (1-r)^\alpha \int_{\mathbb S} |f(rx')|dx' \right)^p \cr
& \times \int_0^1 \left(\int_0^1 \chi_{L_{\epsilon, \alpha}(f)}(r) \frac{(1-r)^{t-\alpha}dr}{(1-r\rho)^{t+1}} \right)^p
(1-\rho)^{\alpha p -1} d\rho\cr
& \leq C_\epsilon \| f \|_{B^{\infty, 1}_\alpha}^p.
\end{align*}
completes the proof. $\Box$

For $f \in h(\mathbb B)$, $\epsilon > 0$ and $\alpha > 0$ we define
$$\hat L_{\epsilon, \alpha}(f) = \{ r : 0 \leq r < 1, M_\infty (f, r) (1-r)^\alpha \geq \epsilon \}.$$
Using Lemma \ref{Lemma3} one easily checks that $B^{p, \infty}_\alpha \subset A^\infty_\alpha$. The following theorem is similar
to the previous one, the difference is that space $B^{\infty, 1}_\alpha$ is replaced by $A^\infty_\alpha$.

\begin{thm}\label{Thm5}
Let $\alpha > 0$ and $1 \leq p < \infty$. For $f \in A^\infty_\alpha$ we set
$$\tilde s_1(f) = {\rm dist}_{A^\infty_\alpha} (f, B^{p, \infty}_\alpha)$$
and
$$\tilde s_2(f) = \inf \left\{ \epsilon > 0 : \int_0^1 \chi_{\hat L_{\epsilon, \alpha}(f)}(r) (1-r)^{-1} dr < \infty \right\}.$$
Then $\tilde s_1(f) \asymp \tilde s_2(f)$, with constants involved depending only on $\alpha$, $p$ and $n$.
\end{thm}

\textbf{Proof. } We use the same ideas as above. We first prove $\tilde s_2(f) \leq \tilde s_1(f)$. Assume
$\tilde s_1(f) < \tilde s_2(f)$. Then there are $0<\epsilon_1 < \epsilon$ and $f_1 \in B^{p, \infty}_\alpha$ such that
$\| f - f_1 \|_{A^\infty_\alpha} \leq \epsilon_1$ and
\begin{equation}\label{contild}
\int_0^1 \chi_{\hat L_{\epsilon, \alpha}(f)}(r) (1-r)^{-1}dr = \infty.
\end{equation}
Since $\| f - f_1 \|_{A^\infty_\alpha} \leq \epsilon_1$ we have
\begin{eqnarray*}
(1-r)^\alpha |f_1(x)| & = & |(1-r)^\alpha f(x) - (1-r)^\alpha [f(x) - f_1(x)] | \cr
& \geq & (1-r)^\alpha |f(x)| - \epsilon_1
\end{eqnarray*}
for every $x \in \mathbb B$, $|x| = r$. Hence, for $r \in \hat L_{\epsilon, \alpha}(f)$, we have
$$M_\infty (f_1, r)(1-r)^\alpha \geq M_\infty (f, r)(1-r)^\alpha - \epsilon_1 \geq \epsilon - \epsilon_1,$$
which implies
$$\chi_{\hat L_{\epsilon, \alpha}(f)}(r) \leq \frac{1}{(\epsilon - \epsilon_1)^p} M_\infty^p(f_1, r)(1-r)^{\alpha p}.$$
Since $f_1 \in B^{p, \infty}_\alpha$ this easily leads to a contradiction with (\ref{contild}).

Next we prove $\tilde s_1(f) \leq C \tilde s_2(f)$. Let us fix $f \in A^\infty_\alpha$ and $\epsilon > 0$ such that
the integral appearing in the definition of $\tilde s_2(f)$ is finite. We use the same decomposition $f = f_1 + f_2$ that appears in (\ref{dec}), and it suffices to prove the following two estimates:
\begin{equation}\label{fdva}
\| f_2 \|_{A^\infty_\alpha} \leq C\epsilon,
\end{equation}
\begin{equation}\label{fjedan}
\| f_1 \|_{B^{p, \infty}_\alpha} \leq C_{\epsilon, p}  \| f \|_{A^\infty_\alpha}.
\end{equation}
Set $\hat L = \hat L_{\epsilon, \alpha}(f)$. Since $|f(\rho y')| \leq \epsilon (1-\rho)^{-\alpha}$ for $\rho \in I \setminus \hat L$ we have,
using Lemma \ref{Lemma1}, for every $x = rx' \in \mathbb B$
\begin{eqnarray*}
|f_2(x)| & \leq & C \epsilon \int_{I \setminus \hat L} \int_{\mathbb S} (1-\rho^2)^\alpha \left(
\frac{(1-r\rho)^{-\{ \alpha \}}}{|\rho x - y'|^{n + [\alpha]}} + \frac{1}{(1-r\rho)^{1+\alpha}} \right) \rho^{n-1} \frac{d\rho dy'}{(1-\rho)^\alpha} \cr
& \leq & C\epsilon \int_0^1 \frac{d\rho}{(1-r\rho)^{1+\alpha}} + C\epsilon \int_0^1 (1-r\rho)^{-\{\alpha\}}
\int_{\mathbb S} \frac{dy'}{|\rho x - y'|^{n+[\alpha]}} d\rho \cr
& \leq &  C\epsilon (1-r)^{-\alpha},
\end{eqnarray*}
which gives estimate (\ref{fdva}). Since $|f(y)| \leq \| f \|_{A^\infty_\alpha} (1-\rho)^{-\alpha}$, $|y| = \rho$, we
have for any $x = rx' \in \mathbb B$:
\begin{eqnarray*}
|f_1(x)| & \leq & \| f \|_{A^\infty_\alpha} \int_{\hat L} \int_{\mathbb S} (1-\rho^2)^\alpha
|Q_\alpha (x, \rho y')| (1-\rho)^{-\alpha} \rho^{n-1} d\rho dy' \cr
& \leq & C \| f \|_{A^\infty_\alpha} \int_{\hat L} \int_{\mathbb S} |Q_\alpha (x, \rho y')| d\rho dy' \cr
& \leq & C \| f \|_{A^\infty_\alpha} \int_{\hat L} \frac{d\rho}{(1-r\rho)^{1+\alpha}},
\end{eqnarray*}
where we used Lemma \ref{Lemma1} again. We proved that $M_\infty (f_1, r) \leq C \| f \|_{A^\infty_\alpha} \phi(r)$, where
\begin{equation}\label{phi}
\phi(r) = \int_{\hat L} \frac{d\rho}{(1-r\rho)^{1+\alpha}}, \qquad 0 \leq r < 1.
\end{equation}
Desired estimate (\ref{fjedan}) will be established once we prove that
$$\psi(r) = (1-r)^\alpha \phi(r) \in L^p(I, (1-r)^{-1}dr).$$
An application of Fubini's theorem gives
\begin{eqnarray*}
\int_0^1 (1-r)^\alpha \int_{\hat L} \frac{d\rho}{(1-r\rho)^{1+\alpha}} \frac{dr}{1-r} & = & \int_{\hat L}
\int_0^1 \frac{(1-r)^{\alpha -1}dr}{(1-r\rho)^{1+\alpha}} d\rho \cr
& \leq & C_\alpha \int_{\hat L} \frac{d\rho}{1-\rho} < \infty
\end{eqnarray*}
by the condition imposed on $\epsilon$. This proves that $\psi \in L^1(I, (1-r)^{-1}dr)$. Since $(1-r)^\alpha (1-\rho)
< (1-r\rho)^{1+\alpha}$ for $r, \rho \in I$ we have
$$\psi(r) = \int_{\hat L} \frac{(1-r)^\alpha (1-\rho)}{(1-r\rho)^{1+\alpha}} \frac{d\rho}{1-\rho} \leq \int_{\hat L}
\frac{d\rho}{1-\rho} = C,$$
and this proves that $\psi \in L^\infty(I, (1-r)^{-1}dr)$. But then it clearly follows that $\psi \in L^p(I, (1-r)^{-1}dr)$ for all $1 \leq p < \infty$. $\Box$

The above theorem can be extended to cover the case $0<p<1$, this is the content of the next theorem.

\begin{thm}
Let $\alpha > 0$ and $0<p\leq 1$. For $f \in A^\infty_\alpha$ we set
$$\overline s_1(f) = {\rm dist}_{A^\infty_\alpha} (f, B^{p, \infty}_\alpha),$$
$$\overline s_2(f) = \inf \left\{ \epsilon > 0 : \int_0^1 \left( \int_{\hat L_{\epsilon, \alpha}(f)}
\frac{dr}{(1-r\rho)^{\alpha+1}}\right)^p (1-\rho)^{\alpha p -1} d\rho < \infty \right\}.$$
Then $\overline s_1(f) \asymp \overline s_2(f)$, with constants involved depending only on $\alpha$, $p$ and $n$.
\end{thm}

\textbf{Proof. } As usual, we prove $\overline s_2(f) \leq \overline s_1(f)$ arguing by contradiction: if $\overline s_1(f) < \overline s_2(f)$ then there are $0 < \epsilon_1 < \epsilon$ and $f_1 \in B^{p, \infty}_\alpha$ such that
$ \| f - f_1 \|_{A^\infty_\alpha} \leq \epsilon_1$ and
\begin{equation}\label{hatl}
J_\epsilon = \int_0^1 \left( \int_{\hat L_\epsilon} \frac{dr}{(1-r\rho)^{\alpha+1}} \right)^p (1-\rho)^{\alpha p - 1} d\rho = \infty,
\end{equation}
where $\hat L_\epsilon = \hat L_{\epsilon, \alpha}(f)$. Since
\begin{eqnarray*}
(1-r)^\alpha M_\infty (f_1, r) & \geq & (1-r)^\alpha M_\infty (f, r) - (1-r)^\alpha M_\infty (f-f_1, r) \cr
& \geq & (1-r)^\alpha M_\infty(f, r) - \epsilon_1,
\end{eqnarray*}
we have $(1-r)^\alpha M_\infty (f_1, r) \geq \epsilon - \epsilon_1$ for $r \in \hat L_\epsilon$, or
$$\chi_{\hat L_\epsilon}(r) \leq (\epsilon - \epsilon_1)^{-1} (1-r)^\alpha M_\infty (f_1, r).$$
Therefore, using Lemma \ref{Lemma2}, we obtain:
\begin{eqnarray*}
J_\epsilon & \leq & C \int_0^1 \left( \int_0^1
\frac{(1-r)^\alpha M_\infty(f_1, r)}{(1-r\rho)^{\alpha + 1}} dr \right)^p (1-\rho)^{\alpha p -1} d \rho \cr
& \leq & C \int_0^1\int_0^1 \frac{(1-r)^{p(\alpha +1) - 1} M_\infty (f_1, r)^p}{(1-r\rho)^{(\alpha+1)p}} (1-\rho)^{\alpha p - 1} dr d\rho \cr
& \leq & C \int_0^1 M_\infty (f_1, r)^p (1-r)^{\alpha p-1} dr = C \| f_1 \|_{B^{p, \infty}_\alpha} < \infty
\end{eqnarray*}
which contradicts (\ref{hatl}). Next we turn to the estimate $\overline s_1(f) \leq C \overline s_2(f)$, using the same
technique as in the previous theorems: we fix $f \in A^\infty_\alpha$ and choose $\epsilon > 0$ such that the
integral appearing in the definition of $\overline s_2(f)$ is finite. Again we use decomposition $f = f_1 + f_2$ from
(\ref{dec}) and the argument presented in the proof of Theorem \ref{Thm5} gives $\| f_2 \|_{A^\infty_\alpha} \leq C \epsilon$, see derivation of (\ref{fdva}). Now it clearly suffices to prove $f_1 \in B^{p, \infty}_\alpha$, in fact we show that
\begin{equation}\label{fj}
\| f_1 \|_{B^{p, \infty}_\alpha} \leq C_\epsilon \| f \|_{A^\infty_\alpha}.
\end{equation}
As in the proof of the previous theorem we have $M_\infty (f_1, r) \leq C \| f \|_{A^\infty_\alpha} \phi(r)$, where
$\phi$ is defined by (\ref{phi}). Therefore (\ref{fj}) follows from the following
$$
\int_0^1 \phi(r)^p (1-r)^{\alpha p-1} dr  =  \int_0^1 \left( \int_{\hat L_\epsilon}
\frac{d\rho}{(1-r\rho)^{1+\alpha}} \right)^p (1-r)^{\alpha p-1} dr = C_\epsilon < \infty.\quad \Box$$

Let us introduce, for $0 < r <1$, $\alpha > -1$ and $f \in h(\mathbb B)$
$$A_\alpha (f, r) = \int_{|w| \leq r|} |f(w)|(1-|w|)^\alpha dw.$$
Next we consider, for $\alpha > -1$, $\beta > 0$ and $0 < p < \infty$ the spaces
$$M_\beta^\alpha(\mathbb B) = \left\{ f \in h(\mathbb B) : \| f \|_{\alpha, \beta} = \sup_{0 \leq r < 1} (1-r)^\beta
A_\alpha (f, r) < \infty \right\},$$
$$M_{p, \beta}^\alpha (\mathbb B) = \left\{ f \in h(\mathbb B) : \| f \|_{\alpha, \beta; p}^p = \int_0^1
(1-r)^{\beta p -1} A_\alpha (f, r)^p dr < \infty \right\}.$$
The spaces $M_\beta^\alpha$ and $M_{p, \beta}^\alpha$, $1 \leq p < \infty$, are Banach spaces and the spaces
$M_{p, \beta}^\alpha$ are complete metric spaces for $0 < p < 1$. Using Lemma \ref{Lemma3} it is readily seen that
$B_{p, \beta}^\alpha \subset B_\beta^\alpha$. We note that the above spaces were investigated in the case of
dimension two by several authors, for example in \cite{JPS} and \cite{SM}. In order to investigate distance
problem in these spaces we introduce, for $\epsilon > 0$, $\alpha > -1$, $\beta > 0$ and $f \in h(\mathbb B)$ the set
$$\tilde L_{\epsilon, \beta}^\alpha (f) = \left\{ r \in (0, 1) : (1-r)^\beta A_\alpha(f, r) \geq \epsilon \right\}.$$

\begin{thm}
Let $p \geq 1$, $\alpha > -1$, $\beta > 0$. We set, for $f \in M_\beta^\alpha$,
$$ t_1(f) = {\rm dist}_{M_\beta^\alpha} (f, M_{p, \beta}^\alpha),$$
$$t_2(f) = \inf \left\{ \epsilon > 0 : \int_0^1 \chi_{\tilde L_{\epsilon, \beta}^\alpha (f)}(r) (1-r)^{-1}dr < \infty
\right\}.$$
Then $t_1(f) \asymp t_2(f)$, with constants involved depending only on $\alpha$, $\beta$, $p$ and $n$.
\end{thm}

We omit the proof of this theorem since it is quite similar to the proofs of all theorems presented in this section.

\end{document}